\documentclass[12pt]{amsart}%

\usepackage{amsmath}
\usepackage{graphicx}
\usepackage{amscd}
\usepackage{geometry}
\usepackage{amsfonts}
\usepackage{amssymb}%
\usepackage{amsthm,enumerate,graphicx,ifsym,mathrsfs}
\usepackage{hyperref}

\newtheorem{theorem}{Theorem}[section]
\theoremstyle{plain}
\newtheorem{Theorem}{Theorem}[section]	
\theoremstyle{plain}

\newtheorem{Conjecture}{Conjecture}

\newtheorem{Definition}{Definition}
\newtheorem{Example}{Example}

\newtheorem{Lemma}[theorem]{Lemma}

\newcommand{\BN}{\mathbb{N}}
\newcommand{\BZ}{\mathbb{Z}}
\newcommand{\BQ}{\mathbb{Q}}
\newcommand{\BR}{\mathbb{R}}
\newcommand{\BC}{\mathbb{C}}

\newcommand{\BI}{\mathbb{I}}

\newenvironment{Proof}{\par\noindent{\sc Proof}\quad}{\hfill\qed\par\smallskip}
\newenvironment{RestateTheorem}[3]{\par\vspace{12pt}\noindent{\bf #1~\ref{#2}}(#3){\bf .}\it}{\par\vspace{12pt}}
\let\oldtocsection=\tocsection
\let\oldtocsubsection=\tocsubsection
\let\oldtocsubsubsection=\tocsubsubsection
\renewcommand{\tocsection}[2]{\hspace{0em}\oldtocsection{#1}{#2}}
\renewcommand{\tocsubsection}[2]{\hspace{2em}\oldtocsubsection{#1}{#2}}
\renewcommand{\tocsubsubsection}[2]{\hspace{4em}
\oldtocsubsubsection{#1}{#2}}

\begin{document}

\title{A Necessary and Sufficient Condition for a Self-Diffeomorphism of a Smooth Manifold to be the Time-1 Map of the Flow of a Differential Equation}
\author{Jeffrey J. Rolland }
\address{Unaffiliated}
\email{rollandj@uwm.edu}
\date{\today}

\begin{abstract}

In topological dynamics, one considers a topological space $X$ and a self-map $f: X \to X$ of $X$ and studies the self-map's properties. In global analysis and geometric mechanics, one considers a smooth manifold $M^n$ and a differential equation $\xi: M \to TM$ on $M$ and studies the flow $\Phi_t: M \times \BR \to M$ of the differential equation. In this paper, we consider a necessary and sufficient condition for a self-diffeomorphism $f$ of a manifold $M$ to be the time-1 map $\Phi_1$ of the flow of a differential equation on $M$.

\end{abstract}
\maketitle

\section{Introduction and Main Result \label{Section: Introduction and Main Result}}

Let $M^n$ be a smooth manifold, let $h: M \to \BR^n$ be a smooth map, and let $p \in M$ be an isolated zero of $h$. Then we might attempt to use \textbf{\textit{Newton's method}} to find $p$, that is, we might choose $x_0$ in a neighborhood of $p$, find a solution $v_i$ of $h(x_i) + Dh(x_i)(v) = 0$, then set $x_{i+1} = \exp(v_i)$ and hope $\lim\limits_{i \to \infty} (x_i) = p$. (A modified version of this might be used to find an isolated critical point of a differential equation on a Riemannian manifold $M$.) If $Dh(p): T_p(M) \to \BR^n$ is nonsingular for all $p \in M$, define $f: M \to M$ by $f(p) = \exp\{Dh(p)^{-1}[-h(p)]\}$. Then $f$ is a self-map of $M$ of the kind that might be considered in smooth dynamics. One wants to know the behavior of the sequences $(x_i) = f^i(x_0)$ for various choices of $x_0$, in particular, when they converge to such an isolated zero $p$, when they tend towards a ``cyclic sequence'' of points $(p_1, p_2, \ldots, p_k)$, and when and how quickly they diverge to infinity. Various sets, such as the Fatou set and the Julia set of $f$, can then be defined.

Let $M^n$ be a smooth manifold, let $\mathfrak{X}(M) = \{\zeta: M \to TM\ |\ \pi_M \circ \zeta = \text{id}_M\}$ be the set of all the differential equations (also known as a tangent vector fields) on $M$, and let $\xi \in \mathfrak{X}(M)$. (We are most interested in the case that $M = TQ$ for $Q$ a closed, Riemannian manifold -- the configuration space [or ``c-space''] of a robot arm, for instance; see \cite{L-P} -- and $\xi$ is the vector field corresponding to an ``elementary'' Lagrangian on $TQ$ [elementary with respect to a specific coordinate patch on $Q$]; see, for instance, \cite{Arnold}). Then, if $\xi$ is \textbf{\textit{complete}}, $\xi$ has a flow,  $\Phi_t: M \times \BR \to M$; see, for instance, \cite{A-R}. Note 1) that $\xi$ is always complete if $M$ is closed, and 2) $\Phi_1 = f$ is a self-diffeomorphism of $M$.

We are interested in when a self-diffeomorphism $f$ of a manifold $M$ is the time-1 map $\Phi_1$ of the flow $\Phi_t$ of a differential equation $\xi$ on $M$. It is known that the set of such self-diffeomorphisms of a manifold $M$ form a Baire first category set from \cite{Palis}.

Here is the statement of the main results.

\begin{Theorem}[A Necessary and Sufficient Condition for a Self-Diffeomorphism of a Manifold to be the Time-1 Map of the Flow of a Differential Equation] \label{thmflow}
Let $n \in \BZ^+$, $M^n$ be a closed, smooth manifold, and $f$ be a self-diffeomorphism of $M$. Then $f = \Phi_1$ for $\Phi_t$ the flow of a differential equation on $M$ if and only if $f$ is smoothly rootable to the identity. In this case, for a smooth root system to the identity $(g_b)$ and corresponding flow $\Psi_t$, the differential equation $\xi$ inducing $\Psi_t$ is given by $\displaystyle \xi = \frac{\partial}{\partial t}\left\{ \Psi_t \right\}\big|^{t=0}$.
\end{Theorem}

\begin{Conjecture}[Rolland Rigidity] \label{conjrigidity}
Suppose $M$ is a locally symmetric space, $\rho$ is the natural Riemannian metric on $M$, and $\nabla$ is the Levi-Civita connection associated to $\rho$ on $M$. We suspect that there is a class of self-diffeomorphisms $f$ of $M$ with the property that if there are two distinct smooth root systems to the identity $(g_{1,b})$ and $(g_{2,b})$ inducing distinct vector fields $\xi_1$ and $\xi_2$, then for each $p \in M$, $p$ and $f(p)$ are conjugate points of $\nabla$. For this class of self-diffeomorphisms, there should be a kind of \underline{rigidity} of $M$ determined by $f$, that is, there should be a differential equation $\xi_b$ on $M$ (called a {\bf base differential equation of $f$}), positive integers $k_1$ and $k_2$, and inversion symmetries $P_1$ and $P_2$ of $(M,\rho)$ with $k_iDP_i \circ \xi_c = \xi_i \circ P_i$, where $DP_i: TM \to TM$ is the induced map on the tangent bundle. In particular, if $\nabla$ is nonpositively curved, $(M, \nabla)$ has no conjugate points, and for any self-diffeomorphism $f$ of this class, there should be a unique differential equation $\xi$ inducing a necessarily unique flow $\Phi_t$ with $f = \Phi_{t=1}$.
\end{Conjecture}

The author would like to thank Lee Mosher of Rutgers University and Jason DeVito of University of Tennessee at Martin for useful conversations on the website \\ \url{https://math.stackexchange.com}.

\section{Proof of the Main Result \label{Section: Proof of Main Result}}
\begin{Definition}
Let $M^n$ be a smooth manifold. Let $f: M \to M$ be a self-diffeomorphism of $M$. We say $f$ is \textbf{rootable to the identity} if $f$ is isotopic to the identity (note this implies $f$ is orientation-preserving if $M$ is orientable) and there is a sequence of self-diffeomorphisms $(g_b)$ of $M$ with 
\begin{enumerate}[{1)}]
\item each $g_b$ isotopic to the identity (note this implies each $g_b$ is orientation-preserving if $M$ is orientable)
\item $(g_b)^b = f$ (that is, each $g_b$ is a ``$b^{\text{th}}$ root of $f$'')
\item (\textbf{The commutativity condition}) for $b_1, b_2 \in \BN$, $g_{b_1}[g_{b_2}(p)] = g_{\frac{b_1+b_2}{b_1b_2}}(p) = g_{\frac{b_2+b_1}{b_2b_1}}(p) = g_{b_2}[g_{b_1}(p)]$
\item (\textbf{The coherency condition}) for $a \in \BZ$, $(g_b)^a = \left(g_{\frac{b}{GCD(a,b)}}\right)^{\frac{a}{GCD(a,b)}}$
\item $\lim\limits_{b \to \infty} (g_b) = \text{id}_M$.
\end{enumerate}

We call such a sequence $(g_b)$ \textbf{a root system to the identity}. Note one should only need to determine $(g_b)$ on some cofinal subset of the naturals that leads to a dense subset of the rationals, e.g., $b = 2^c$, leading to the dyadic rationals. 

Set $\Phi_{\frac{a}{b}}: M \times \BQ \to M$ by $\Phi_{\frac{a}{b}}(p) = (g_b)^a(p)$. We say $f$ is \textbf{continuously rootable to the identity} if $\Phi_{\frac{a}{b}}$ extends to to a continuous function $\Phi_{t}: M \times \BR \to M$ and call the root system to the identity $(g_b)$ a \textbf{continuous root system to the identity} in this case. Note that if such a continuous extension $\Phi_t$ exists, it is unique. We say $f$ is \textbf{smoothly rootable to the identity} if $\Phi_t$ is smooth and call the root system to the identity $(g_b)$ a \textbf{smooth root system to the identity} in this case. 
\end{Definition}

\begin{Lemma} \label{commutePower}
Let $(g_b)$  be a root system to the identity, let $a_1, a_2 \in \BZ$, and let $b_1, b_2 \in \BZ^+$. Then $g_{b_2}^{a_2}[g_{b_1}^{a_1}(p)] = (g_{b_1b_12})^{a_1b_2+a_1b_2}(p)]$
\end{Lemma}

\begin{Proof}
$\displaystyle g_{b_2}^{a_2}[g_{b_1}^{a_1}(p)] = \displaystyle g_{b_1b_2}^{a_2b_1}[g_{b_1b_2}^{a_1b_2}(p)]$ 
by the coherency condition, so
$\displaystyle g_{b_2}^{a_2}[g_{b_1}^{a_1}(p)]=\displaystyle (g_{b_1b_2})^{a_2b_1+a_1b_2}(p)$ 
by the definition of composition of functions.\\
\end{Proof}

\begin{RestateTheorem}{Theorem}{thmflow}{A Necessary and Sufficient Condition for a Self-Diffeomorphism of a Manifold to be the Time-1 Map of the Flow of a Differential Equation} 
Let $n \in \BZ^+$, $M^n$ be a closed, smooth manifold, and $f$ be a self-diffeomorphism of $M$. Then $f = \Phi_1$ for $\Phi_t$ the flow of a differential equation on $M$ if and only if $f$ is smoothly rootable to the identity. In this case, for a smooth root system to the identity $(g_b)$ and corresponding flow $\Psi_t$, the differential equation $\xi$ inducing $\Psi_t$ is given by $\displaystyle \xi = \frac{\partial}{\partial t}\left\{ \Psi_t \right\}\big|^{t=0}$.
\end{RestateTheorem}

\begin{Proof} Let $n \in \BZ^+$, $M^n$ be a closed, smooth manifold, and $f$ be a self-diffeomorphism of $M$. \\
($\Rightarrow$) Suppose $f = \Phi_1$ for $\Phi_t$ the flow of a differential equation $\xi$ on $M$. Define $H: M \times \BI \to M$ by $H(p,t) = \Phi(p,1-t)$. Then $H(p,0) = \Phi_{1-0}(p) = \Phi_{1}(p) = f(p)$ and $H(p,1) = \Phi_{1-1}(p) = \Phi_{0}(p) = \text{id}_M(p)$, so $f$ is isotopic to $\text{id}_M$. Let $b \in \BZ^+$ and set $g_b(p) = \Phi_{\frac{1}{b}}(p)$.

\begin{enumerate}[{1)}]
\item Define $H_b: M \times \BI \to M$ by $H_b(p,t) = \Phi(p,\frac{1-t}{b})$. Then $H(p,0) = \Phi_{\frac{1-0}{b}}(p) = \Phi_{\frac{1}{b}}(p) = g_b(p)$ and $H_b(p,1) = \Phi_{\frac{1-1}{b}}(p) = \Phi_{0}(p) = \text{id}_M(p)$, so each $g_b$ is isotopic to $\text{id}_M$.

\item $(g_b)^b(p) = \Phi_{\frac{1}{b}}^b(p) = \Phi_{b\frac{1}{b}}(p) = \Phi_1(p) = f(p)$, so each $g_b$ is a ``$b^{\text{th}}$ root of $f$''.

\item Let $b_1, b_2 \in \BZ^+$. Then $g_{b_1}[(g_{b_2})(p)] = \Phi_{\frac{1}{b_1}}[\Phi_{\frac{1}{b_2}}(p)] = \Phi_{\frac{1}{b_1}+\frac{1}{b_2}}(p)] \text{ as } \Phi_t \text{ is a flow } = \Phi_{\frac{1}{b_2}+\frac{1}{b_1}}(p) = \Phi_{\frac{1}{b_2}}[\Phi_{\frac{1}{b_1}}(p)] = g_{b_2}[(g_{b_1})(p)]$, so $(g_b)$ obeys the commutativity condition.

\item For $a \in \BZ$, $(g_b)^a(p) = (\Phi_{\frac{1}{b}})^a(p) = \Phi_{\frac{a}{b}}(p) = \Phi_{\frac{\frac{a}{GCD(a,b)}}{\frac{b}{GCD(a,b)}}}(p) = \left(g_{\frac{b}{GCD(a,b)}}\right)^{\frac{a}{GCD(a,b)}}(p)$ as $\Phi_t$ is a flow, so $(g_b)$ obeys the coherency condition.

\item $\lim\limits_{b \to \infty}(g_b)(p) = \lim\limits_{b \to \infty}\Phi_{\frac{1}{b}}(p) = \Phi_0(p) = \text{id}_M$.
\end{enumerate}

\noindent
Hence, $(g_b)$ is a smooth root system to the identity, and $f$ is smoothly rootable to the identity.

($\Leftarrow$) Suppose $f$ is smoothly rootable to the identity. Let $(g_b)$ be a smooth root system to the identity for $f$. Define $\Psi_{\frac{a}{b}}: M \times \BQ \to M$ by $\Psi_{\frac{a}{b}}(p) = (g_b)^a(p)$. Then $\left(g_{\frac{b}{GCD(a,b)}}\right)^{\frac{a}{GCD(a,b)}}(p) = \Psi_{\frac{\frac{a}{GCD(a,b)}}{\frac{b}{GCD(a,b)}}}(p)$ by the coherency condition, so $\Psi_{\frac{a}{b}}$ is well-defined. Moreover, as $(g_b)$ is a smooth root system to the identity, it is a continuous root system to the identity, so $\Psi_{\frac{a}{b}}$ extends uniquely to a continuous function $\Psi_t: M \times \BR \to M$. Next, as $(g_b)$ is a smooth root system to the identity, $\Psi_t: M \times \BR \to M$ is smooth. We must show $\Psi_t$ is a flow on $M$. \\

\begin{enumerate}[{a)}]
\item Note $\Psi_0 = \Psi_{\frac{0}{b}} = (g_{b})^0 = \text{id}_M$. 

\item
\begin{enumerate}[{(i)}]
\item Let $a_1, a_2 \in \BZ$ and $b_1, b_2 \in \BZ^+$. Then 

\begin{tabular}{ll}
$\Psi_{\frac{a_2}{b_2}}\left[\Psi_{\frac{a_1}{b_1}}(p)\right]$ & $=(g_{b_1})^{a_1}[(g_{b_2})^{a_2}(p)]$ \\
 & $=(g_{b_1b_2})^{a_1b_2+a_2b_1}(p)$ by Lemma \ref{commutePower}  \\
 & $=$ $\Psi_{\frac{a_1b_2+a_2b_1}{b_1b_2}}(p)$ \\
 & $=$ $\Psi_{\frac{a_1b_2}{b_1b_2}+\frac{a_2b_1}{b_1b_2}}(p)$ \\
$\Psi_{\frac{a_2}{b_2}}\left[\Psi_{\frac{a_1}{b_1}}(p)\right]$ & $=\Psi_{\frac{a_1}{b_1}+\frac{a_2}{b_2}}(p) $ by the fact that $\Psi_{\frac{a}{b}}$ is well-defined on $\BQ$. 
\end{tabular}

\item It then follows that $\Psi_{t_2}\left[\Psi_{t_1}(p)\right] = \Psi_{t_1+t_2}(p)$ for $t_1, t_2 \in \BR$ by continuity as $(g_b)$ is a continuous root system to the identity.

\item The fact that $\Psi_t$ is smooth  follows from the fact that $\Psi_{t}$ is a smooth root system to the identity.
\end{enumerate}
\end{enumerate}

\noindent
This shows that $\Psi_t$ is a flow and $f = \Psi_1$.

A proof of the fact that for a complete, time-independent vector field $\xi$, the differential equation $\xi$ inducing $\Psi_t$ is given by $\displaystyle \xi = \frac{\partial}{\partial t}\left\{ \Psi_t \right\}\big|^{t=0}$ may be found in, for instance, \cite{Ar-Go}.
\end{Proof}

\section{Examples \label{Section: Examples}}
\begin{Example}
Let $M = S^1$ and let $f$ be the antipodal map, $f(p) = -p = e^{i\pi}p = e^{-i\pi}p$. Then there are two obvious smooth root systems to the identity, $\displaystyle g_{1,b}(p) = e^{\frac{i\pi}{b}}p$ and $\displaystyle g_{2,b}(p) = e^{\frac{-i\pi}{b}}p$. Note that the isometry $P$ of $S^1$ with its standard Riemannian metric defined by $P(p) = \overline{p}$ (complex conjugation, picturing $S^1$ as being the unit circle in $\BC$) has $DP \circ\xi_1 = \xi_2 \circ P$, for $\xi_i$ the differential equation determined by $(g_{i,b})$, where $DP: TS^1 \to TS^1$ is the induced map on the tangent bundle.
\end{Example}

\begin{Example}
(A generalization of an example from Jason Devito) With $S^3$, thinking of $S^3$ as a Lie group, the antipodal map (left multiplication by -1, $L_{-1}$) has uncountably many square roots: left multiplication by any purely imaginary unit quaternion. Every imaginary quaternion has exactly two quaternion square roots, $u_3$ and $-u_3$ with $(\pm u_3)^2 = q$. Only one of the $\pm u_3$'s at level $3$ will have a angle smaller that $q$ with 1, the other one will be $-u_3$ and will have a smaller angle than $q$ with -1. This pattern continues with $u_{c-1}$  has exactly two quaternion square roots, $u_c$ and $-u_c$ with $(\pm u_c)^2 = u_{c-1}$. Only one of the $\pm u_c$'s at level $c$ will have a angle smaller that $u_{c-1}$ with 1, the other one will be $-u_c$ and will have a smaller angle than $-u_{c-1}$ with -1. If $g_{2^c} = L_{u_c}$, then $(g_{2^c})$ is a sequence of $2^{c \text{ th}}$ roots of $f = L_{-1}$ defined on a cofinal subset of the naturals with each $(g_{2^c})$ a smooth root system to the identity. Hence, we have a case where we have uncountably many different differential equations $\xi_q$ with $\Phi_{q,t=1} = f$. Note that if $q_1$ and $q_2$ are distinct imaginary quaternions with $q_1^2 = q_2^2 = -1$, the isometry $P$ of $S^3$ defined by $P(p) = L_{q_2(q_1^{-1})}(p)$has $DP \circ \xi_1 = \xi_2 \circ P$, for $\xi_i$ the differential equation determined by $(g_{i,b})$, where $DP: TS^3 \to TS^3$ is the induced map on the tangent bundle.
\end{Example}

\begin{RestateTheorem}{Conjecture}{conjrigidity}{Rolland Rigidity}
Suppose $M$ is a locally symmetric space, $\rho$ is the natural Riemannian metric on $M$, and $\nabla$ is the Levi-Civita connection associated to $\rho$ on $M$. We suspect that there is a class of self-diffeomorphisms $f$ of $M$ with the property that if there are two distinct smooth root systems to the identity $(g_{1,b})$ and $(g_{2,b})$ inducing distinct vector fields $\xi_1$ and $\xi_2$, then for each $p \in M$, $p$ and $f(p)$ are conjugate points of $\nabla$. For this class of self-diffeomorphisms, there should be a kind of \underline{rigidity} of $M$ determined by $f$, that is, there should be a differential equation $\xi_b$ on $M$ (called a {\bf base differential equation of $f$}), positive integers $k_1$ and $k_2$, and inversion symmetries $P_1$ and $P_2$ of $(M,\rho)$ with $k_iDP_i \circ \xi_c = \xi_i \circ P_i$, where $DP_i: TM \to TM$ is the induced map on the tangent bundle. In particular, if $\nabla$ is nonpositively curved, $(M, \nabla)$ has no conjugate points, and for any self-diffeomorphism $f$ of this class, there should be a unique differential equation $\xi$ inducing a necessarily unique flow $\Phi_t$ with $f = \Phi_{t=1}$.
\end{RestateTheorem}

Note that the set $G_f = \{P_i\ |\ 1DP_i \circ \xi_b = \xi_i \circ P_i\}$, for $\xi_i$ a differential equation on $M$ inducing $f$ with coefficient 1 for its inversion symmetry conjugating it to the base differential equation $\xi_b$, is a group under composition of inversion symmetries. For any choice of basepoint $*$ of $M$, $G_f$ embeds in $GL[T_*(M)]$. A choice of base differential equation $\xi_b$ from any of the $\xi_i$'s with $1DP_i \circ \xi_b = \xi_i \circ P_i$ for $P_i \in G_f$ used to define $G_f$ should be thought of as a choice similar to the choice of basepoint in defining the fundamental group.


\begin{thebibliography}{01}
\bibitem{A-R} R. Abraham and J. Robbin. \textit{Transversal mappings and flows}. W. A. Benjamin, 1967.
\bibitem{Ar-Go} J. Arango and A. G\'{o}mez. Flows and diffeomorphisms. Revista Colombiana de Matem´aticas,
32(1):13–27, Jan. 1998.
\bibitem{Arnold} V. Arnold. \textit{Mathematical methods of classical mechanics}, volume 60 of Graduate Texts in Mathematics.
Springer, 1989.
\bibitem{L-P} K. Lynch and F. Park. \textit{Modern robotics: mechanics, planning, and control}. Cambridge University
Press, 1st edition, 2017.
\bibitem{Palis} J. Palis. Vector fields generate few diffeomorphisms. Bulletin of the American Mathematical
Society, 80:503–505, 1974.
\end{thebibliography}
\end{document}